\documentclass[12pt]{amsart}
\usepackage{amssymb}
\usepackage{xcolor}

\usepackage[T1]{fontenc}
\usepackage[cp1250]{inputenc}
\usepackage{amssymb}
\usepackage{amsmath}
\usepackage{amsfonts}
\usepackage{amsthm}
\usepackage{graphics}

\newtheorem{conjecture}{Conjecture}
\newtheorem{theorem}{Theorem}

\begin{document}

\title[Primes in numerical semigroups]{Primes in numerical semigroups}

\thanks{$^1$ Partially supported by Program MATH AmSud, Grant MATHAMSUD 18-MATH-01, Project FLaNASAGraTA.}
\author{J.L. Ram\'irez Alfons\'in$^1$
\and  
M. Ska\l ba }

\address{Institute of Mathematics\\
University of Warsaw\\
Banacha 2\\
02-097 Warszawa, Poland}
\email{skalba@mimuw.edu.pl}

\address{IMAG, Univ.\ Montpellier, CNRS, Montpellier, France and 
UMI2924 - Jean-Christophe Yoccoz, CNRS-IMPA}
\email{jorge.ramirez-alfonsin@umontpellier.fr}



\begin{abstract}
Let $0<a<b$ be two relatively prime integers and let $\langle a,b\rangle$ be the numerical semigroup generated by $a$ and $b$ with {\em Frobenius number} $g(a,b)=ab-a-b$. In this note, we prove that there exists a prime number $p\in \langle a,b\rangle$ with $p<g(a,b)$ when the product $ab$ is sufficiently large. Two related conjectures are posed and discussed as well. 
\end{abstract}

\maketitle

Let $0<a<b$ be two relatively prime integers. Let $S=\langle a,b\rangle=\{n \ | \ n=ax+by, x,y\in\mathbb{Z}, x,y\ge 0\}$ be the numerical semigroup generated by $a$ and $b$. A well-known result due to Sylvester \cite{sylv} states that the largest integer not belonging to $S$, denoted by $g(a,b)$, is given by $ab-a-b$. $g(a,b)$ is called the {\em Frobenius number}. We refer the reader to \cite{ramirez1} for a complexity result on the Frobenius number for general numerical semigroups and to \cite{ramirez2} for an extensive literature on it. 

We clearly have that any prime $p$ larger than $g(a,b)$ belongs to $\langle a,b\rangle$. A less obvious and more intriguing question is whether there is a prime $p \le g(a,b)$ belonging to $\langle a,b\rangle$.  

In this note, we show that there always exists a prime $p\in \langle a,b\rangle, p <g(a,b)$ when the product $ab$ is sufficiently large. The latter is a straight forward consequence of the below Theorem.

Let $0<u<v$ be integers. We define 
$$\pi_S[u,v] =|\{p \ \text{prime}\ |\ p \in S, \ u\le p\le v\}|.$$ 
For short, we may write $\pi_S$ instead of $\pi_S[0,g(a,b)]$.
 
\begin{theorem}\label{theo:main} Let $3\le a<b$ be two relatively prime integers and let $S=\langle a,b\rangle$
be the numerical semigroup generated by $a$ and $b$. Then, for any fixed $\varepsilon>0$ there exists $C(\varepsilon)>0$ such that
$$\pi_S >C(\varepsilon) \frac{g(a,b)}{\log(g(a,b))^{2+\varepsilon}}$$
for $ab$ sufficiently large.
\end{theorem}

Let us quickly introduce some notation and recall some facts needed for the proof of Theorem \ref{theo:main}.

Let $S=\langle a,b\rangle$ and let $0<u<v$ be integers. We define 
$$n_S[u,v]= |\{n\in \mathbb{N}\ |\ u\le n\le v, \ n\in S\}|$$ 
and
$$n^c_{S}[u,v]= |\{n\in\mathbb{N} \ |\ u\le n\le v, \ n\not\in S\}|.$$ 

For short, we may write $n_S$ instead of $n_S[0,g(a,b)]$ and $n^c_S$ instead of $n_S^c[0,g(a,b)]$.
The set of elements in $n_S^c=\mathbb{N}\setminus S$ are usually called the {\em gaps} of $S$.
\smallskip

It is known \cite{ramirez2}  that $S$ is always {\em symmetric}, that is, for any integer $0\le s \le g(a,b)$
$$\hbox{$s\in S$ if and only if $g(a,b)-s\not\in S$}.$$
Obtaining that 
$$n_S=\frac{g(a,b)+1}{2}.$$

Let $\pi(x)$ be the number of primes integers less or equals to $x$ for any real number $x$. We recall that the well-known prime number theorem asserts that
\begin{equation}\label{eq:prime}
\pi(x)\sim \frac{x}{\log x}.
\end{equation}
\smallskip

We may now prove Theorem \ref{theo:main}.
\smallskip

{\em Proof of Theorem \ref{theo:main}.} Let $\varepsilon>0$ be fixed. We distinguish two cases.

{\bf Case 1)} Suppose that $a> (\log(ab))^{1+\varepsilon}$. Let us take $c=ab/(\log(ab))^{1+\varepsilon}$. It is known \cite{MRT} that if $k\in [0,\dots ,g(a,b)]$ then 

$$n_S[0,k]=\sum\limits_{i=0}^{\lfloor\frac{k}{b}\rfloor}\left(\left\lfloor\frac{k-ib}{a}\right\rfloor+1\right).$$ 

In our case, we obtain that

$$\begin{array}{ll}
n_S[0,c] & <\lfloor\frac{c}{a}\rfloor+\lfloor\frac{c}{b}\rfloor\left(\lfloor\frac{c-b}{a}\rfloor+1\right)
 < \lfloor\frac{c}{a}\rfloor+\lfloor\frac{c}{b}\rfloor\left(\lfloor\frac{c}{a}\rfloor+1\right)\\
 \\
 & < \frac{c}{a}+\frac{c}{b}+ \frac{c^2}{ab}< \frac{bc+ac+c^2}{ab}< \frac{2c^2+c^2}{ab}=\frac{3c^2}{ab}=\frac{3ab}{(\log(ab))^{2+2\varepsilon}}
\end{array}$$
 where the last inequality holds since $c>b>a$.

Due to the symmetry of $S$, we have

\begin{equation}\label{eq:nonrep}
n_S^c[g(a,b)-c,g(a,b)]=n_S[0,c]<\frac{3ab}{(\log(ab))^{2+2\varepsilon}}.
\end{equation}

Now, by \eqref{eq:prime}, we have

\begin{equation}\label{eq:difprimes}
\pi_S[g(a,b)-c,g(a,b)]=\pi(g(a,b))-\pi(g(a,b)-c)>>\frac{c}{\log (ab)}=\frac{ab}{(\log(ab))^{2+\varepsilon}}
\end{equation}
when $ab$ is large enough.

Finally, by combining equations \eqref{eq:nonrep} and \eqref{eq:difprimes}, we obtain

$$\begin{array}{ll}
\pi_S& > \pi_S[g(a,b)-c,g(a,b)]-n^c_S[g(a,b)-c,g(a,b)]\\
& >\frac{ab}{(\log(ab))^{2+\varepsilon}}-\frac{3ab}{(\log(ab))^{2+2\varepsilon}}>0
\end{array}$$

where the last inequality holds since $(\log(ab))^{\varepsilon}>3$ for $ab$ large enough for the fixed $\epsilon$. The above leads to the desired estimation of $\pi_S$.
\medskip

{\bf Case 2)} Suppose that $3\leq a\leq  (\log(ab))^{1+\varepsilon}$.

If $p\in [b,\dots ,g(a,b)]$ is a prime and $p\equiv b\pmod{a}$ then $p$ is clearly representable as $p=b+\frac{p-b}{a} a$. By Siegel-Walfisz theorem \cite{prachar,walf}, the number of such primes $p$, denoted by $N$, is
$$ N=\frac{1}{\varphi(a)}\int_b^{g(a,b)}\frac{du}{\log u}+R $$
where $\varphi$ is the {\em Euler totient function} and $|R|< D'(\varepsilon)\frac{g(a,b)}{(\log(g(a,b)))^{2+2\varepsilon}}$ uniformly in $[a,\dots ,g(a,b)]$. 

Since the function  $1/\log u$ is decreasing on the interval $[b,g(a,b)]$ then
$$\int_b^{g(a,b)}\frac{du}{\log u}>(g(a,b)-b)\cdot \frac{1}{\log g(a,b)}$$
and therefore

\begin{equation}\label{eq:N}
N> \frac{1}{\varphi(a)}\cdot \frac{g(a,b)-b}{\log(g(a,b))}-D'(\varepsilon)\frac{g(a,b)}{(\log(g(a,b)))^{2+2\varepsilon}}.
\end{equation}


Now, we have that 
$$\begin{array}{l}
\frac{1}{\varphi(a)}\cdot \frac{g(a,b)-b}{\log(g(a,b))} \cdot
\frac{(\log(g(a,b)))^{2+\varepsilon}}{g(a,b)}\\
=\frac{1}{\varphi(a)}\log(g(a,b))^{1+\varepsilon}\left( 1-\frac{b}{g(a,b)}\right)\\
>\frac{1}{\log(ab)^{1+\varepsilon}}\log(g(a,b))^{1+\varepsilon}\left( 1-\frac{b}{g(a,b)}\right) \hfill \hbox{(since $(\log(ab))^{1+\varepsilon}\ge a >\varphi(a)$)} \\
>\left( \frac{\log(ab)-\log(3)}{\log(ab)}\right)^{1+\varepsilon} \frac{1}{5}>F>0\hfill \hbox{(since $g(a,b))> ab/3$ and $\frac{b}{g(a,b)}\le \frac{4}{5}$)}
\end{array}$$

for some absolute $F>0$, uniformly for $ab\geq D''(\varepsilon)$ with $a\geq 3$. 

Yielding to 
\begin{equation}\label{eq:N1}
\frac{1}{\varphi(a)}\cdot \frac{g(a,b)-b}{\log(g(a,b))}  \geq 
F \frac{g(a,b)}{\log(g(a,b))^{2+\varepsilon}}
\end{equation}

and combining equations \eqref{eq:N} and \eqref{eq:N1} we obtain 
$$N>F'\frac{g(a,b)}{\log(g(a,b))^{2+\varepsilon}} $$
for $ab$ large enough for the fixed $\epsilon$. The latter leads to the desired estimation of $\pi_S$ also in this case.

\hfill\ensuremath{\blacksquare}
 \smallskip

\section{Concluding remarks}

A number of computer experiments lead us to the following. 

\begin{conjecture}\label{c1} Let $1\le a<b$ be two relatively prime integers and let $S$ be the numerical semigroup generated by $a$ and $b$. Then,
$$\pi_S>0.$$
\end{conjecture}

{In analogy with the symmetry of $\langle a, b\rangle$ mentioned above, our task of looking for primes in $\langle a, b\rangle$ is related with the task of finding primes in $[g(a,b)-1)/2,\ldots, g(a,b)]$. From this point of view, Conjecture \ref{c1} can be thought of as a counterpart of the famous Chebyshev theorem stating that there is always a prime in $[n, \dots, ,2n]$\ for any $n\geq 2$, see \cite[Chapter 3]{sierp}. 
A way to attack Conjecture \ref{c1} could be by applying {\em effective versions} of Siegel-Walfisz theorem. For instance, one may try to use \cite[Corollary 8.31]{tenen} in order to get computable constants in our estimates. However, it is not an easy task to trace all constants appearing in the relevant estimates of $L(x,\chi)$ (but in principle possible). The remaining cases for {\em small} values $ab$ must to be treated by computer. 

\begin{conjecture} Let $1\le a<b$ be two relatively prime integers and let $S$ be the numerical semigroup generated by $a$ and $b$. Then,
$$\pi_S\sim \frac{\pi(g(a,b))}{2}\ \mbox{ for }ab\rightarrow\infty. $$
\end{conjecture}

In the same spirit as the prime number theorem, this conjecture seems to be out of reach. 
\smallskip

The famous Linnik's theorem 
asserts that there exist absolute constants $C$ and $L$ such that: for given relatively prime integers $a,b$ the least prime $p$ satisfying $p\equiv b\pmod{a}$ is less than $Ca^L$.
It is conjectured that the value $L=2$, but the current record is only that $L\leq 5$, see  \cite{xylou}.
\smallskip


On the same flavor of Linnik's theorem that concerns the existence of primes of the form $ax+b$, Theorem \ref{theo:main} is concerning the existence of primes of the form $ax+by$ with $x,y\geq 1$ less than $ab$ for sufficiently large $ab$. This relation could shed light on in either direction.

\end{document}